 \documentclass{monsky2009} 

\usepackage{graphicx}

\usepackage{amssymb}
\usepackage{amsmath}

\usepackage{bm}

\newenvironment{theorem*}[1]{\textbf{#1}\itshape \hspace{.3em}}{\upshape}
\newenvironment{remark*}[1]{\textbf{#1}\itshape \hspace{.3em}}{\upshape}
\newenvironment{corollary*}[1]{\textbf{#1}\itshape \hspace{.3em}}{\upshape}
\newenvironment{proof}{\textbf{Proof\hspace{.3em}}}{}

\newtheorem{definition}{Definition}[section]
\newtheorem{theorem}[definition]{Theorem}
\newtheorem{lemma}[definition]{Lemma}

\newtheorem{conjecture}[definition]{Conjecture}

\newcommand{\charstic}{\ensuremath{\mathrm{char\ }}}

\newcommand{\bilprod}{\mbox{\ensuremath\,\#\,}}

\begin{document}

\begin{frontmatter}






\title{Algebraicity of some Hilbert-Kunz multiplicities (modulo a conjecture)}
\author{Paul Monsky}

\address{Brandeis University, Waltham MA  02454-9110, USA. monsky@brandeis.edu}

\begin{abstract}
Let $F$ be a finite field of characteristic 2 and $h$ be the element $x^{3}+y^{3}+xyz$ of $F[[x,y,z]]$. In an earlier paper we made a precise conjecture as to the values of the colengths of the ideals $(x^{q},y^{q},z^{q},h^{j})$ for $q$ a power of 2. We also showed that if the conjecture holds then the Hilbert-Kunz series of $H=uv+h$ is algebraic (of degree 2) over $Q(w)$, and that $\mu (h)$ is algebraic (explicitly, $\frac{4}{3}+\frac{5}{14\sqrt{7}}$). In this note, assuming the same conjecture, we use a theory of infinite matrices to rederive this result, and we extend it to a wider class of $H$; for example $H=g(u,v)+h$. In a follow-up paper, under the same hypothesis, we will show that transcendental Hilbert-Kunz multiplicities exist.

\end{abstract}


\end{frontmatter}


\section{A product on $X$}
\label{section1}

In this section we develop some general results about Hilbert-Kunz series and multiplicities for characteristic 2 power series. (There are similar results, implicit in \cite{5}, in all finite characteristics but they are harder to prove.)

\begin{definition}
\label{def1.1}
$X$ is the vector space of functions $I\rightarrow Q$ where $I=[0,1]\cap Z[\frac{1}{2}]$. If $f\ne 0$ is in the maximal ideal of $F[[u_{1},\ldots , u_{r}]]$, $\charstic F=2$, then $\phi_{f}$ in $X$ is the function $\frac{i}{q}\rightarrow q^{-r}\deg (u_{1}^{q}, \ldots , u_{r}^{q},f^{i})$; here $q$ denotes a power of 2 and $deg$ is colength in $F[[u_{1},\ldots , u_{r}]]$. Note that $\phi_{f}$ is well-defined.
\end{definition}

\begin{definition}
\label{def1.2}
$\alpha$ in $X$ is convex if for all $i$ and $q$ with $0<i<q$, $2\alpha\left(\frac{i}{q}\right)\ge\alpha\left(\frac{i-1}{q}\right)+\alpha\left(\frac{i+1}{q}\right)$.
\end{definition}

Note that $\phi_{f}(0)=0$, $\phi_{f}(1)=1$, $\phi_{f}$ is convex and $\phi_{f}$ is Lipschitz. The first two assertions are clear. If we set $J=(u_{1}^{q}, \ldots , u_{r}^{q})$ then multiplication by $f$ induces a map of $(J,f^{i-1})/(J,f^{i})$ onto $(J,f^{i})/(J,f^{i+1})$, yielding convexity. Finally, as Lipschitz constant we can take the Hilbert-Kunz multiplicity, $\mu$, of~$f$.

\begin{definition}
\label{def1.3}
Suppose that $\alpha$ in $X$ is convex Lipschitz with $\alpha(0)=0$ and $\alpha(1)=1$. Then $\mu(\alpha)=\lim_{n\rightarrow\infty}\alpha(2^{-n})\cdot 2^{n}$, while $\mathcal{S}_{\alpha}$ is the element\linebreak $\sum\alpha(2^{-n})(2w)^{n}$ of $Q[[w]]$. (The convexity of $\alpha$ shows that $n\rightarrow 2^{n}\alpha(2^{-n})$ is non-decreasing. Since $\alpha$ is Lipschitz, the function is bounded and the limit exists.)
\end{definition}

\begin{remark*}{Remarks}
When $\alpha = \phi_{f}$, $\mu(\alpha)$ and $\mathcal{S}_{\alpha}(2^{r-1}w)$ are just the Hilbert-Kunz multiplicity and Hilbert-Kunz series of $f$. Note that if $\alpha$ is as in Definition \ref{def1.3} then $\mu(\alpha)=\lim_{w\rightarrow 1^{-}}(1-w)\mathcal{S}_{\alpha}(w)$. For convexity shows that the co-efficients of the power series $(1-w)\mathcal{S}_{\alpha}(w)$ are $\ge 0$. So the limit is the value of this power series at 1. And we note that $\alpha(1)+\left(2\alpha\left(\frac{1}{2}\right)-\alpha(1)\right) + \left(4\alpha\left(\frac{1}{4}\right)-2\alpha\left(\frac{1}{2}\right)\right) + \cdots$ converges to $\mu(\alpha)$.
\end{remark*}

We next define a bilinear product $\bilprod : X\times X \rightarrow X$ and show that if $f\ne 0$ and $g\ne 0$ are in the maximal ideals of $F[[u_{1},\ldots , u_{r}]]$ and $F[[v_{1},\ldots , v_{s}]]$, then $\phi_{f}\bilprod \phi_{g}=\phi_{h}$, where $h$ is the element $f(u)+g(v)$ of $F[[u_{1},\ldots , u_{r}, v_{1},\ldots , v_{s}]]$. (There is a similar construction, implicit in \cite{5}, in any finite characteristic.)

\begin{definition}
\label{def1.4}
Suppose $\alpha$ and $\beta$ are in $X$.  We define $\alpha \bilprod \beta(t)$ by induction on the denominator of $t$ in $I$, according to the following procedure:

Let $\alpha_{0}$ and $\alpha_{1}$ be the elements $t\rightarrow\alpha\left(\frac{t}{2}\right)$ and $t\rightarrow\alpha\left(\frac{1+t}{2}\right)$ of $X$; define $\beta_{0}$ and $\beta_{1}$ similarly. Then:
\begin{enumerate}
\item[(1)] $\alpha \bilprod \beta(0)=0$\hfill $\alpha \bilprod \beta(1)=(\alpha(1)-\alpha(0))(\beta(1)-\beta(0))$\hfill\rule{0pt}{0pt}
\item[(2)] If $0\le t \le \frac{1}{2}$

$\alpha\bilprod\beta(t)=\alpha_{0}\bilprod\beta_{0}(2t)+\alpha_{1}\bilprod\beta_{1}(2t)$
\item[(3)] If $\frac{1}{2}\le t \le 1$

$\alpha\bilprod\beta(t)=\alpha_{0}\bilprod\beta_{0}(1)+\alpha_{1}\bilprod\beta_{1}(1)+\alpha_{0}\bilprod\beta_{1}(2t-1)+\alpha_{1}\bilprod\beta_{0}(2t-1)$
\end{enumerate}
\end{definition}

Note that when $t=0$, $\frac{1}{2}$ or $1$ the two definitions of $\alpha\bilprod\beta(t)$ given by the above scheme coincide, so that $\alpha\bilprod\beta$ is a well-defined element of $X$. $\bilprod$ is evidently bilinear and symmetric; one can show that it is associative. It's easy to see that if $\alpha$ is constant then $\alpha\bilprod\beta=0$, while if $\alpha$ is the identity function $t$, $\alpha\bilprod\beta=(\beta(1)-\beta(0))t$. In particular, $t\bilprod t=t$.

Now let $T_{0}$ and $T_{1}$ $X\rightarrow X$ be the maps taking $\alpha$ to $t\rightarrow \alpha\left(\frac{t}{2}\right)$ and $t\rightarrow \alpha\left(\frac{1+t}{2}\right)$. Replacing $t$ by $\frac{t}{2}$ in (2) above and by $\frac{1+t}{2}$ in (3) above gives:

\begin{theorem}
\label{theorem1.5}
If $\gamma=\alpha\bilprod\beta$ then:
\begin{enumerate}
\item[]$T_{0}(\gamma)=(T_{0}(\alpha)\bilprod T_{0}(\beta))+(T_{1}(\alpha)\bilprod T_{1}(\beta))$
\item[]$T_{1}(\gamma)=\gamma\left(\frac{1}{2}\right)+(T_{0}(\alpha)\bilprod T_{1}(\beta))+(T_{1}(\alpha)\bilprod T_{0}(\beta))$
\end{enumerate}
\end{theorem}

We now recall some notation used in both \cite{1} and \cite{5}. By an $F[T]$-module we'll mean a finitely generated $F[T]$-module annihilated by a power of $T$. $\Gamma$ is the Grothendieck group of the set of isomorphism classes of such modules. There is a multiplication on $\Gamma$ making it into a commutative ring; if $V$ and $W$ are $F[T]$-modules, a representative of their product is $V\underset{F}{\otimes}W$, with $T$ acting by $(T_{V}\otimes \mathrm{id})+(\mathrm{id}\otimes T_{W})$. There is a $Z$-basis $\lambda_{0}, \lambda_{1}, \ldots$ of $\Gamma$ with the following property. If $V$ is an $F[T]$-module then the class of $V$ in $\Gamma$ is $\sum c_{i}\lambda_{i}$ where $c_{i}=(-1)^{i}\dim(T^{i}V/T^{i+1}V)$.  Because $\charstic F=2$, the multiplicative structure of $\Gamma$ is very simple; $\lambda_{i}\lambda_{j}=\lambda_{k}$ where $k$ is the ``Nim-sum'' of $i$ and $j$.

\begin{definition}
\label{def1.6}
If $\alpha$ is in $X$, $n\ge 0$ and $q=2^{n}$, then $\mathcal{L}_{n}(\alpha)$ is the element $\sum_{0}^{q-1}\left(\alpha\left(\frac{i+1}{q}\right)-\alpha\left(\frac{i}{q}\right)\right)(-)^{i}\lambda_{i}$ of $\Gamma_{Q}= \Gamma\underset{Z}{\otimes}Q$.
\end{definition}

Note that $\mathcal{L}_{0}(\alpha)=(\alpha(1)-\alpha(0))\lambda_{0}$. If $\alpha=\phi_{f}$ and $V$ is the $F[T]$-module $F[[u_{1}^{q},\ldots , u_{r}^{q}]]$ with $T$ acting by multiplication by $f$, then $q^{r}\mathcal{L}_{n}(\alpha)=\linebreak \sum \dim (T^{i}V/T^{i+1}V)(-)^{i}\lambda_{i}$; this is precisely the class of $V$ in $\Gamma$.

Suppose now that $q=2^{n}$ and $0\le i < q$. Since the Nim-sum of $i$ and $q$ is $q+i$, $\lambda_{i}\lambda_{q}=\lambda_{q+i}$ giving:

\begin{lemma}
\label{lemma1.7}
$\mathcal{L}_{n+1}(\alpha)=\mathcal{L}_{n}(\alpha_{0})+\lambda_{q}\mathcal{L}_{n}(\alpha_{1})$
\end{lemma}

\begin{theorem}
\label{theorem1.8}
If $\gamma=\alpha\bilprod\beta$, $\mathcal{L}_{n}(\gamma)=\mathcal{L}_{n}(\alpha)\cdot\mathcal{L}_{n}(\beta)$.
\end{theorem}

\begin{proof}
We argue by induction on $n$. Since $\gamma(1)-\gamma(0)=(\alpha(1)-\alpha(0))(\beta(1)-\beta(0))$ the result holds for $n=0$. Suppose that it's true for a given $n$. Lemma \ref{lemma1.7}, Theorem \ref{theorem1.5} and the induction hypothesis show that $\mathcal{L}_{n+1}(\gamma)=\linebreak\mathcal{L}_{n}(\alpha_{0})\mathcal{L}_{n}(\beta_{0})+\mathcal{L}_{n}(\alpha_{1})\mathcal{L}_{n}(\beta_{1})+\lambda_{q}(\mathcal{L}_{n}(\alpha_{0})\mathcal{L}_{n}(\beta_{1})+\mathcal{L}_{n}(\alpha_{1})\mathcal{L}_{n}(\beta_{0}))$.  But this is $(\mathcal{L}_{n}(\alpha_{0})+\lambda_{q}\mathcal{L}_{n}(\alpha_{1}))\cdot(\mathcal{L}_{n}(\beta_{0})+\lambda_{q}\mathcal{L}_{n}(\beta_{1}))$ which is $\mathcal{L}_{n+1}(\alpha)\cdot\mathcal{L}_{n+1}(\beta)$ by Lemma \ref{lemma1.7}.
\qed
\end{proof}

\begin{theorem}
\label{theorem1.9}
Suppose $h=f(u_{1},\ldots , u_{r})+g(v_{1},\ldots , v_{s})$. Then:
\begin{enumerate}
\item[(1)] For each $n$, $\mathcal{L}_{n}(\phi_{h})=\mathcal{L}_{n}(\phi_{f})\cdot\mathcal{L}_{n}(\phi_{g})$
\item[(2)] $\phi_{h} = \phi_{f}\bilprod\phi_{g}$
\end{enumerate}
\end{theorem}

\begin{proof}
With $q=2^{n}$, let $V$ be as in the paragraph following Definition \ref{def1.6}. As we've seen $V$ represents the element $q^{r}\mathcal{L}_{n}(\phi_{f})$ of $\Gamma$. Replacing $f$ by $g$ we get a $W$ representing the element $q^{s}\mathcal{L}_{n}(\phi_{g})$ of $\Gamma$. Then $q^{r+s}\mathcal{L}_{n}(\phi_{f})\cdot\mathcal{L}_{n}(\phi_{g})$ is represented by $F[[u_{1},\ldots , u_{r},v_{1},\ldots ,v_{s}]]/(u_{1}^{q},\ldots ,v_{s}^{q})$ with $T$ acting by multiplication by $f(u_{1},\ldots , u_{r})+g(v_{1},\ldots , v_{s})=h$. Since this $F[T]$-module represents $q^{r+s}\mathcal{L}_{n}(\phi_{h})$ we get (1). Suppose now that $\phi_{h}(t)\ne \phi_{f}\bilprod\phi_{g}(t)$ for some $t=\frac{i}{2^{n}}$. Choose such a $t$ with $i$ as small as possible. Then $i\ne 0$, and the co-efficients of $\lambda_{i-1}$ in $\mathcal{L}_{n}(\phi_{h})$ and $\mathcal{L}_{n}(\phi_{f}\bilprod\phi_{g})$ differ. Theorem \ref{theorem1.8} then shows that $\mathcal{L}_{n}(\phi_{h})\ne \mathcal{L}_{n}(\phi_{f})\mathcal{L}_{n}(\phi_{g})$, contradicting (1).
\qed
\end{proof}
\pagebreak

\begin{theorem}
\label{theorem1.10}
If $\alpha$ and $\beta$ are Lipschitz with Lipschitz constant $m$, then $\gamma=\alpha\bilprod\beta$ is Lipschitz with Lipschitz constant $m^{2}$.
\end{theorem}

\begin{proof}
We show that if $0\le j <2q$ then $\left|\gamma\left(\frac{j+1}{2q}\right)-\gamma\left(\frac{j}{2q}\right)\right|\le \frac{m^{2}}{2q}$, arguing by induction on $q$. Note first that $\alpha_{0}$, $\alpha_{1}$, $\beta_{0}$ and $\beta_{1}$ are all Lipschitz with Lipschitz constant $\frac{m}{2}$. We claim that when $j< q$ the values of $\alpha_{0}\bilprod\beta_{0}$ (and of $\alpha_{1}\bilprod\beta_{1}$) at $\frac{j+1}{q}$ and $\frac{j}{q}$ differ by at most $\frac{m^{2}}{4q}$. (When $q=1$, $j=0$, and this is clear. When $q>1$ we use the fact that $\alpha_{0}$ and $\beta_{0}$ (and $\alpha_{1}$ and $\beta_{1}$) have Lipschitz constant $\frac{m}{2}$, together with the induction hypothesis.)  Theorem \ref{theorem1.5} then shows that $\gamma\left(\frac{j+1}{2q}\right)$ and $\gamma\left(\frac{j}{2q}\right)$ differ by at most $\frac{m^{2}}{4q}+\frac{m^{2}}{4q}=\frac{m^{2}}{2q}$. The argument is similar when $j\ge q$, but now we make use of the values of $\alpha_{0}\bilprod\beta_{1}$ (and of $\alpha_{1}\bilprod\beta_{0}$) at $\frac{j+1-q}{q}$ and $\frac{j-q}{q}$.
\qed
\end{proof}

\begin{lemma}
\label{lemma1.11}
Let $\delta_{r}$, $r\ge 1$, be the class of $F[T]/T^{r}$ in $\Gamma$; note that $\delta_{r}=\lambda_{0}-\lambda_{1}+\lambda_{2}\cdots +(-)^{r-1}\lambda_{r-1}$. Then for $\alpha$ in $X$ the following are equivalent:
\begin{enumerate}
\item[(1)] $\alpha$ is convex.
\item[(2)] For each $n$, $\mathcal{L}_{n}(\alpha)=\sum_{0}^{q-1}c_{i}(-)^{i}\lambda_{i}$ with $c_{0}\ge c_{1}\ge \cdots \ge c_{q-1}$.
\item[(3)] For each $n$, $\mathcal{L}_{n}(\alpha)$ is a linear combination of $\delta_{1},\ldots , \delta_{q}$ with the co-efficients of $\delta_{1}, \ldots , \delta_{q-1} \ge 0$.
\end{enumerate}
\end{lemma}

\begin{proof}
Since the $c_{i}$ in (2) is $\alpha\left(\frac{i+1}{q}\right)-\alpha\left(\frac{i}{q}\right)$, (1) and (2) are equivalent. Suppose (2) holds. If we set $c_{q}=0$, then the formula for $\delta_{r}$ given above shows that $\mathcal{L}_{n}(\alpha)=\sum_{0}^{q-1}(c_{i}-c_{i+1})\delta_{i}$. Since $c_{0}-c_{1}, \cdots , c_{q-2}-c_{q-1}$ are all $\ge 0$ we get (3). That (2) follows from (3) is easy.
\qed
\end{proof}

\begin{lemma}
\label{lemma1.12}
Suppose $1\le r,s\le q$. Then, in $\Gamma$, $\delta_{r}\delta_{s}$ is a linear combination of $\delta_{1},\ldots , \delta_{q}$ with non-negative integer co-efficients. Furthermore $\delta_{r}\delta_{q}=r\delta_{q}$.
\end{lemma}

\begin{proof}
Let $V$ and $W$ be the $F[T]$-modules $F[T]/T^{r}$ and $F[T]/T^{s}$ representing $\delta_{r}$ and $\delta_{s}$. Writing  $V\underset{F}{\otimes}W$ (with $T$ acting by $T_{V}\otimes\mathrm{id} + \mathrm{id}\otimes T_{W}$) as a direct sum of cyclic $F[T]$-modules we get the first assertion. The second is an easy calculation.
\qed
\end{proof}

\begin{theorem}
\label{theorem1.13}
If $\alpha$ and $\beta$ in $X$ are convex, then so is $\alpha\bilprod\beta$.
\end{theorem}

\begin{proof}
By Lemma \ref{lemma1.11}, $\mathcal{L}_{n}(\alpha)$ and $\mathcal{L}_{n}(\beta)$ are each linear combinations of $\delta_{1},\ldots , \delta_{q}$ with the co-efficients of $\delta_{1},\ldots , \delta_{q-1}\ge 0$. By Lemma \ref{lemma1.12} the same is true of $\mathcal{L}_{n}(\alpha)\cdot\mathcal{L}_{n}(\beta)$. Theorem \ref{theorem1.8} and Lemma \ref{lemma1.11} then show that $\alpha\bilprod\beta$ is convex.
\qed
\end{proof}

\begin{theorem}
\label{theorem1.14}
Suppose that $\alpha$ in $X$ is convex Lipschitz with $\alpha(0)=0$ and $\alpha(1)=1$. Suppose further that $\mathcal{S}_{\alpha}=\sum\alpha(2^{-n})(2w)^{n}$ lies in a finite extension, $L$, of $Q(w)$. (We extend the imbedding of $Q[w]$ in $Q[[w]]$ to their fields of fractions.)  Then $\mu(\alpha)$ is algebraic over $Q$ of degree $\le [L: Q(w)]$. In fact there is a valuation ring containing $Q[w]$ in $L$ whose maximal ideal contains $w-1$ and whose residue class field contains a copy of $Q(\mu(\alpha))$.
\end{theorem}

\begin{proof}
Take $H$ irreducible in $Q[W,T]$ so that $H(w,(1-w)\mathcal{S}_{\alpha})=0$. Then for any $z$ in the open unit disc, $H(z,(1-z)\mathcal{S}_{\alpha}(z))=0$. The remarks following Definition \ref{def1.3} show that $H(1,\mu(\alpha))=0$. Since $H(1,T)\ne 0$, $\mu(\alpha)$ is algebraic over $Q$. Let $g$ be $\mathrm{Irr}(\mu(\alpha),Q)$. Then $(W-1,g(T))$ is a maximal ideal in $Q[W,T]/H$ and we take a valuation ring in $L$ that contains $Q[w,(1-w)\mathcal{S}_{\alpha}]=Q[W,T]/H$, and whose maximal ideal contracts to the above maximal ideal.
\qed
\end{proof}

\section{A calculation from \cite{2}, revisited}
\label{section2}

Let $f$ be the element $x^{3}+y^{3}+xyz$ of $Z/2\,[x,y,z]$, defining a nodal cubic. The values of $\phi_{f}$ at $\frac{1}{q}$ are known, and in particular, $\mu(f)=\frac{7}{3}$. In \cite{2} we conjectured a precise value for all $\phi_{f}\left(\frac{i}{q}\right)$, and showed that the conjecture implied that $\mu(uv+f)$ is $\frac{4}{3}+\frac{5}{14\sqrt{7}}$. In this section we'll rework this result using infinite matrix techniques from \cite{3}; this approach will give rise to more general theorems.

\begin{definition}
\label{def2.1}
$1$, $t$ and $\epsilon$ will denote the elements\, $t\rightarrow 1$, $t\rightarrow t$ and $t\rightarrow t-t^{2}$ of $X$.
\end{definition}

\begin{definition}
\label{def2.2}
For $m=0,1,2,\ldots $ and $t$ in $I$, $\phi_{m}(t)$ is defined by induction on the denominator of $t$ as follows:
\begin{enumerate}
\item[(1)] $\phi_{m}(0)=\phi_{m}(1)=0$
\item[(2)] If $0\le t \le \frac{1}{2}$, $8\phi_{m}(t)=\phi_{m+1}(2t)+(8m+6)t$ for $m$ even, and $\phi_{m-1}(2t)+\epsilon(2t)+(8m+6)t$ for $m$ odd.
\item[(3)] If $\frac{1}{2}\le t \le 1$, $8\phi_{0}(t)=\phi_{0}(2t-1)+6(1-t)$
\item[(4)] If $\frac{1}{2}\le t \le 1$, $8\phi_{m}(t)=\phi_{m-1}(2t-1)+\epsilon(2t-1)+(8m+6)(1-t)$ for $m\ne 0$ even, and $\phi_{m+1}(2t-1)+(8m+6)(1-t)$ for $m$ odd.
\end{enumerate}
\end{definition}

When $t=0$, $\frac{1}{2}$ or $1$, the two definitions of $\phi_{m}(t)$ given by the above scheme evidently coincide. So the $\phi_{m}$ are well-defined elements of $X$. Replacing $t$ by $\frac{t}{2}$ in (2) and by $\frac{1+t}{2}$ in (3) and (4) we get the ``magnification rules'':
\begin{enumerate}
\item[(1)] \parbox{2.5in}{$8T_{0}(\phi_{0})=\phi_{1}+3t$}\hfill
\parbox{2.5in}{$8T_{1}(\phi_{0})=\phi_{0}+3(1-t)$}

\item[(2)] When $m\ne 0$ is even,\\
\parbox{2.5in}{$8T_{0}(\phi_{m})=\phi_{m+1}+(4m+3)t$}\hfill
\parbox{2.5in}{$8T_{1}(\phi_{m})=\phi_{m-1}+\epsilon +(4m+3)(1-t)$}

\item[(3)] When $m$ is odd,\\
\parbox{2.5in}{$8T_{0}(\phi_{m})=\phi_{m-1}+\epsilon +(4m+3)t$}\hfill
\parbox{2.5in}{$8T_{1}(\phi_{m})=\phi_{m+1}+(4m+3)(1-t)$}

\end{enumerate}

Note also that $4T_{0}(\epsilon)=\epsilon +t$ and that $4T_{1}(\epsilon)=\epsilon +(1-t)$.

\begin{conjecture}
\label{conjecture2.3}
If $f=x^{3}+y^{3}+xyz$, then $\phi_{f}=t+\phi_{0}$ with $\phi_{0}$ as above.
\end{conjecture}

In \cite{2} we presented evidence for a conjecture easily seen to be equivalent to this. We noted in particular that both sides agree at all $\frac{1}{q}$ and at each $\frac{i}{512}$.

\begin{theorem}
\label{theorem2.4}
If $E_{1}=\epsilon\bilprod\phi_{0}$ then $\lim_{n\rightarrow\infty}E_{1}(2^{-n})2^{n}=\frac{1}{3}+\frac{5}{14\sqrt{7}}$.
\end{theorem}

Suppose now that Conjecture \ref{conjecture2.3} holds. Then $t+E_{1}=(t+\epsilon)\bilprod(t+\phi_{0})=\phi_{uv}\bilprod\phi_{f}=\phi_{uv+f}$. So Theorem \ref{theorem2.4} tells us that the Hilbert-Kunz multiplicity of $uv+x^{3}+y^{3}+xyz$ is $\lim_{n\rightarrow\infty}\left(2^{-n}+E_{1}(2^{-n})\right)2^{n}=\frac{4}{3}+\frac{5}{14\sqrt{7}}$, an observation made in \cite{2}. We now give a proof of Theorem \ref{theorem2.4} using the techniques of \cite{3}.

\begin{lemma}
\label{lemma2.5}
Let $T:X\rightarrow X$ be $32T_{0}$. Set $E_{k}=\epsilon\bilprod\phi_{k-1}$. Then:

\begin{enumerate}
\item[(1)] $T(E_{1})=E_{1}+E_{2}+6t$
\item[(2)] $T(E_{k})=E_{k-1}+E_{k+1}+(8k-2)t+(\epsilon\bilprod\epsilon)$ for $k>1$
\item[(3)] $T(\epsilon\bilprod\epsilon)=4(\epsilon\bilprod\epsilon)+4t$, and $T(t)=16t$
\end{enumerate}
\end{lemma}

\begin{proof}
Suppose $k$ is even. Then $T(E_{k})=32T_{0}(\epsilon\bilprod\phi_{k-1})=(4T_{0}(\epsilon)\bilprod \linebreak 8T_{0}(\phi_{k-1}))+(4T_{1}(\epsilon)\bilprod 8T_{1}(\phi_{k-1}))$. The magnification rules following Definition \ref{def2.2} show that this is $(\epsilon + t)\bilprod (\phi_{k-2}+\epsilon + (4k-1)t)+(\epsilon + 1-t)\bilprod (\phi_{k}+(4k-1)(1-t))$. Expanding out we get $(\epsilon \bilprod \phi_{k-2})+(4k-1)t+(\epsilon\bilprod\phi_{k})+(4k-1)t+(\epsilon\bilprod\epsilon )=E_{k-1}+E_{k+1}+(8k-2)t+(\epsilon\bilprod\epsilon)$. The other parts of the lemma are derived similarly.
\qed
\end{proof}

\begin{lemma}
\label{lemma2.6}
Let $\mathcal{S}$ be the power series $\sum E_{1}(2^{-n})(32w)^{n}$. Then $(1-16w)(1-4w)(1-2w)^{2}\mathcal{S}=4w(1-2w)^{2}+(2w-12w^{2})\sqrt{1-4w^{2}}$.
\end{lemma}

\begin{proof}
Let $l:X\rightarrow Q$ be evaluation at $1$, so that $l(E_{k})=0$ for each $k$, and $l(\epsilon\bilprod\epsilon)=0$, while $l(t)=1$. Then $E_{1}(2^{-n})32^{n}$ is $l(T^{n}(E_{1}))$ and $\mathcal{S}$ is just $\sum l(T^{n}(E_{1}))w^{n}$. If we take $Y$ to be the subspace of $X$ spanned by $\epsilon\bilprod\epsilon$ and $t$, Lemma \ref{lemma2.5} shows that we are in the situation of Example 5.12 of \cite{3}. The final line of that paper is the desired result.
\qed
\end{proof}

Theorem \ref{theorem2.4} is now easily proved. Lemma \ref{lemma2.6} shows that the value, $\lambda$, of $(1-16w)\mathcal{S}$ at $w=\frac{1}{16}$ is $\left(\frac{4}{3}\cdot\frac{64}{49}\right)\left(\frac{4}{16}\cdot\frac{49}{64}+\frac{5}{64}\sqrt{\frac{63}{64}}\right)=\frac{1}{3}+\frac{5}{14\sqrt{7}}$. Furthermore, $\mathcal{S}-\frac{\lambda}{1-16w}$ is holomorphic in the disc $|w|<\frac{1}{4}$. It follows that $\mathcal{S}\left(\frac{w}{16}\right)-\frac{\lambda}{1-w}$ is holomorphic in $|w|<4$, and so the co-efficients in its power series expansion $\rightarrow 0$. So $E_{1}(2^{-n})\cdot 2^{n}-\lambda \rightarrow 0$, the desired result.
\qed

We conclude this section by showing that the $\phi_{m}$ of Definition \ref{def2.2} are convex and Lipschitz.

\begin{lemma}
\label{lemma2.7}
$\phi_{m}\left(\frac{1}{q}\right)\le \frac{4m+4}{3q}$ for even $m$ and $\frac{4m+3}{3q}$ for odd $m$.
\end{lemma}

\begin{proof}
When $q=2$, $\phi_{m}\left(\frac{1}{q}\right)=\frac{4m+3}{4q}$. We argue by induction. Suppose $q\ge 2$. If $m$ is even, $\phi_{m}\left(\frac{1}{2q}\right)=\frac{1}{8}\phi_{m+1}\left(\frac{1}{q}\right)+\frac{4m+3}{8q}$. By the induction hypothesis this is $\le \frac{4m+7}{24q}+\frac{4m+3}{8q}=\frac{4m+4}{3(2q)}$. If $m$ is odd, $\phi_{m}\left(\frac{1}{2q}\right)=\frac{1}{8}\phi_{m+1}\left(\frac{1}{2q}\right)+\frac{1}{8q}-\frac{1}{8q^{2}}+\frac{4m+3}{8q}$. By the induction hypothesis this is $\le \frac{4m}{24q}+\frac{4m+4}{8q}=\frac{4m+3}{3(2q)}$.
\qed
\end{proof}

\begin{lemma}
\label{lemma2.8}
$\phi_{m}\left(1-\frac{1}{q}\right)\le \frac{4m+4}{3q}$ for odd $m$ and $\frac{4m+3}{3q}$ for even $m$.
\end{lemma}

\begin{proof}
$q=2$ is clear. Suppose $q\ge 2$; we argue by induction. If $m$ is odd, $\phi_{m}\left(1-\frac{1}{2q}\right)=\frac{1}{8}\phi_{m+1}\left(1-\frac{1}{q}\right)+\frac{4m+3}{8q}$, while if $m\ne 0$ is even, $\phi_{m}\left(1-\frac{1}{2q}\right)=\frac{1}{8}\phi_{m+1}\left(1-\frac{1}{q}\right)+\frac{1}{8q}-\frac{1}{8q^{2}}+\frac{4m+3}{8q}$, and we continue as in the proof of Lemma \ref{lemma2.7}. Finally, $\phi_{0}\left(1-\frac{1}{2q}\right) = \frac{1}{8}\phi_{0}\left(1-\frac{1}{q}\right)+\frac{3}{8q}$. By the induction hypothesis this is $\le \frac{1}{8q}+\frac{3}{8q}=\frac{1}{2q}$.
\qed
\end{proof}

\begin{lemma}
\label{lemma2.9}
$\phi_{m}\left(\frac{q+1}{2q}\right)$ and $\phi_{m}\left(\frac{q-1}{2q}\right)$ are $\le\phi_{m}\left(\frac{1}{2}\right)$.
\end{lemma}

\begin{proof}
If $m$ is odd, $8\left(\phi_{m}\left(\frac{1}{2}\right)-\phi_{m}\left(\frac{q+1}{2q}\right)\right)=\frac{4m+3}{q}-\phi_{m+1}\left(\frac{1}{q}\right)$.  By Lemma \ref{lemma2.7} this is $\ge \frac{4m+3}{q}-\frac{4m+8}{3q}\ge 0$. Also $8\left(\phi_{m}\left(\frac{1}{2}\right)-\phi_{m}\left(\frac{q-1}{2q}\right)\right)\ge \frac{4m+3}{q}-\phi_{m-1}\left(1-\frac{1}{q}\right)$. By Lemma \ref{lemma2.8} this is $\ge \frac{4m+3}{q}-\frac{4m-1}{3q}\ge 0$. The argument for even $m$ is similar.
\qed
\end{proof}

\begin{theorem}
\label{theorem2.10}
The $\phi_{m}$ are convex and Lipschitz.
\end{theorem}

\begin{proof}
To prove convexity, we show that if $0< j < 2q$, then $2\phi_{m}\left(\frac{j}{2q}\right)-\phi_{m}\left(\frac{j-1}{2q}\right)-\phi_{m}\left(\frac{j+1}{2q}\right)\ge 0$, arguing by induction on $q$. The case $q=1$ is immediate. When $j<q$ the induction assumption tells us that $2\phi_{s}\left(\frac{j}{q}\right)-\phi_{s}\left(\frac{j-1}{q}\right)-\phi_{s}\left(\frac{j+1}{q}\right)\ge 0$ for each $s$; this and the fact that $\epsilon$ and $t$ are convex gives the result. When $j>q$, the induction assumption tells us that $2\phi_{s}\left(\frac{j-q}{q}\right)-\phi_{s}\left(\frac{j-1-q}{q}\right)-\phi_{s}\left(\frac{j+1-q}{q}\right)\ge 0$; this and the convexity of $\epsilon$ and $1-t$ give the result. Finally the case $j=q$ is handled by Lemma \ref{lemma2.9}. Note also that Lemmas \ref{lemma2.7} and \ref{lemma2.8} show that $|\phi_{m}\left(\frac{1}{q}\right)-\phi_{m}(0)|$ and $|\phi_{m}\left(1-\frac{1}{q}\right)-\phi_{m}(1)|$ are each $\le \frac{4m+4}{3q}$. Since $\phi_{m}$ is convex, it follows that it is Lipschitz with Lipschitz constant $\frac{4m+4}{3}$. 
\qed
\end{proof}

\section{Algebraicity results}
\label{section3}

We generalize the calculations of Section \ref{section2} to show:

\begin{theorem}
\label{theorem3.1}
Suppose $\beta_{1}$ lies in a finite dimensional subspace of $X$ stable under $T_{0}$ and $T_{1}$, and is convex Lipschitz. Set $E_{1}=\beta_{1}\bilprod\phi_{0}$ with $\phi_{0}$ as in Definition \ref{def2.2}. Then the power series $\mathcal{S}_{t+E_{1}}(w)$ is algebraic over $Q(w)$, and $\mu(t+E_{1})$ is algebraic over $Q$.
\end{theorem}

\begin{proof}
Since $\beta_{1}$ and $\phi_{0}$ are convex Lipschitz, the same is true of $t+E_{1}$. In view of Theorem \ref{theorem1.14} we only need to prove the result for $\mathcal{S}$. We shall mimic the proof of Theorem \ref{theorem2.4}. Take $\beta_{1},\ldots \beta_{l},1,t$ spanning a space stable under $T_{0}$ and $T_{1}$. We are free to modify each $\beta_{j}$ by a linear combination of $1$ and $t$ and so may assume $\beta_{j}(0)=\beta_{j}(1)=0$. Then $T_{0}(\beta_{j})=$ (a linear combination of $\beta_{i}$) $+$ a multiple of $t$, while $T_{1}(\beta_{j})=$ (a linear combination of $\beta_{i}$) $+$ a multiple of $(1-t)$. Since $T_{0}(\beta_{j})(1)=T_{1}(\beta_{j})(0)=\beta_{j}\left(\frac{1}{2}\right)$ we get:
\begin{eqnarray*}
T_{0}(\beta_{j}) &=& \sum r_{i,j}\beta_{i}+c_{j}t\\
T_{1}(\beta_{j}) &=& \sum s_{i,j}\beta_{i}+c_{j}(1-t)
\end{eqnarray*}
with the $r_{i,j}$, the $s_{i,j}$ and the $c_{j}$ all in $Q$.

We proceed in several steps:
\begin{enumerate}
\item[I)] Let $R$ and $S$ be the elements $|r_{i,j}|$ and $|s_{i,j}|$ of $M_{l}(Q)$. We define an infinite matrix $V$ with rows and columns indexed by the positive integers as follows. $V$ is built up out of $l$ by $l$ blocks. The initial diagonal block is $S$ while all succeeding diagonal blocks are matrices of zeroes. The blocks just below the diagonal blocks are alternately $R$ and $S$, as are the blocks just to the right of the diagonal blocks. All other entries are zero.

\item[II)] Let $\phi_{m}$ be as in Definition \ref{def2.2}. If $m\ge 0$ and $1 \le j \le l$ let $E_{j+lm}=\beta_{j}\bilprod\phi_{m}$; note that $E_{1}=\beta_{1}\bilprod\phi_{0}$ in accord with the statement of the theorem. $Y\subset X$ is the subspace spanned by $t$ and the $\beta_{j}\bilprod\epsilon$, and we define $y_{1}, y_{2},\ldots $ in $Y$ as follows. If $1\le j \le l$, $y_{j}=6c_{j}t$. If $m>0$, $y_{j+lm}-(8m+6)c_{j}t=\sum r_{i,j}(\beta_{i}\bilprod\epsilon)$ for odd $m$ and $\sum s_{i,j}(\beta_{i}\bilprod\epsilon)$ for even $m$. Note that $4T_{0}(\beta_{j}\bilprod\epsilon)=T_{0}(\beta_{j})\bilprod(\epsilon+t)+T_{1}(\beta_{j})\bilprod(\epsilon+1-t)$, so that $Y$ is stable under $T_{0}$.
\item[III)] With notation as above we claim that $8T_{0}(E_{j})=\sum v_{i,j}E_{i}+y_{j}$. This amounts to:
\begin{enumerate}
\item[(1)] If $1\le j\le l$, $8T_{0}(E_{j})=\sum s_{i,j}E_{i}+\sum r_{i,j}E_{i+l}+y_{j}$
\item[(2)] If $m$ is odd, $8T_{0}(E_{j+lm})=\sum r_{i,j}E_{i+lm-l}+\sum s_{i,j}E_{i+lm+l}+y_{j+lm}$
\item[(3)] If $m>0$ is even, $8T_{0}(E_{j+lm})=\sum s_{i,j}E_{i+lm-l}+\sum r_{i,j}E_{i+lm+l}+y_{j+lm}$
\end{enumerate}
Note that the left hand side of (3) is $8T_{0}(\beta_{j}\bilprod\phi_{m})=(\sum r_{i,j}\beta_{i}+c_{j}t)\bilprod\linebreak(\phi_{m+1}+(4m+3)t)+(\sum s_{i,j}\beta_{i}+c_{j}(1-t))\bilprod(\phi_{m-1}+\epsilon + (4m+3)(1-t))$. Expanding out and using the definition of $y_{j+lm}$ we get (3). Similar calculations give (1) and (2).
\item[IV)] Now set $s=2l$. It's convenient to view the matrix $V$ of I as built up out of $s$ by $s$ blocks. Set $D=\left(\begin{smallmatrix}
S&R\\
R&0
\end{smallmatrix}
\right)$ and $B=\left(\begin{smallmatrix}
0&R\\
R&0
\end{smallmatrix}
\right)$ in $M_{s}(Q)$. Then the diagonal blocks of $V$ are a single $D$ followed by $B$'s. If we take $A=\left(\begin{smallmatrix}
0&S\\
0&0
\end{smallmatrix}
\right)$ and $C=\left(\begin{smallmatrix}
0&0\\
S&0
\end{smallmatrix}
\right)$, then the blocks just below the diagonal blocks are all $A$'s, while those just to the right of the diagonal blocks are all $C$'s. And all other entries are zero.
\end{enumerate}

The proof of Theorem \ref{theorem3.1} is now easy. III and IV tell us that we are in the situation of Theorem 5.11 of \cite{3} with $T=8T_{0}$ and $s$, $A$, $B$, $C$, $D$ as above. (Note that the $y_{j}$ are all in $Y$, that $Y$ is finite-dimensional and stable under $T$, and that the condition of Lemma 5.10 of \cite{3} on the sequence $y_{1},y_{2},\ldots$ is trivially satisfied.) Let $l:X\rightarrow Q$ be evaluation at $1$ so that each $l(E_{j})=0$. Then Theorem 5.11 of \cite{3} shows that $\sum l(T^{n}(E_{1}))w^{n}=\sum E_{1}(2^{-n})(8w)^{n}$ is algebraic over $Q(w)$. So the same is true of  $\frac{1}{1-w}+\sum E_{1}(2^{-n})(2w)^{n}=\sum (2^{-n}+E_{1}(2^{-n}))(2w)^{n}=\mathcal{S}_{t+E_{1}}(w)$.
\qed
\end{proof}

\begin{definition}
\label{def3.2}
$g\ne 0$ in the maximal ideal of $F[[u_{1},\ldots , u_{r}]]$ is ``strongly rational'' if $\phi_{g}$ lies in a finite dimensional subspace of $X$ stable under $T_{0}$ and $T_{1}$.
\end{definition}

 The following is shown in \cite{4} and \cite{5}:
 
 \begin{theorem}
 \label{theorem3.3}
 \rule{0pt}{0pt} 
 \begin{enumerate}
 \item[(1)] If $F$ is finite and $r=2$, $g$ is strongly rational.
 \item[(2)] If $g$ is strongly rational, the Hilbert-Kunz series of $g$ lies in $Q(w)$, and $\mu(g)$ is rational.
 \item[(3)] If $g(u_{1},\ldots , u_{r})$ and $h(v_{1},\ldots , v_{s})$ are strongly rational, then so are $g(u)+h(v)$, $g(u)h(v)$, and all powers of $g(u)$.
 \end{enumerate}
 \end{theorem}

\begin{remark*}{Remark}
Much of the above is easy to prove. (1) however makes use of a result on the finiteness of the  number of ideal classes in certain $1$-dimensional rings. And the proof of (3) for $g(u)+h(v)$ (or rather the generalization of this result to arbitrary finite characteristic $p$) isn't easy. But when $p=2$ there's an immediate proof. Namely suppose that $V_{1}$ and $V_{2}$ are finite dimensional subspaces of $X$ containing $\phi_{g}$ and $\phi_{h}$ and stable under $T_{0}$ and $T_{1}$. Then the space spanned by $1$ and $V_{1}\bilprod V_{2}$ is finite dimensional and stable under $T_{0}$ and $T_{1}$. Furthermore it contains $\phi_{g}\bilprod\phi_{h}=\phi_{g(u)+h(v)}$.
\end{remark*}

If $g$ is strongly rational, Theorem \ref{theorem3.1} tells us that $\mathcal{S}_{t+(\phi_{g}\bilprod\phi_{0})}$ is algebraic over $Q(w)$ and that $\mu(t+(\phi_{g}\bilprod\phi_{0}))$ is algebraic. Now $t+(\phi_{g}\bilprod\phi_{0})=\phi_{g}\bilprod(t+\phi_{0})$. This gives:

\begin{theorem}
\label{theorem3.4}
Suppose that Conjecture \ref{conjecture2.3} holds; that is to say that $t+\phi_{0}=\phi_{x^{3}+y^{3}+xyz}$.  Then if $g$ in $F[[u_{1},\ldots , u_{r}]]$ is strongly rational, the Hilbert-Kunz series of $g(u_{1},\ldots , u_{r})+x^{3}+y^{3}+xyz$ is algebraic over $Q(w)$, and the Hilbert-Kunz multiplicity is algebraic. In particular using Theorem \ref{theorem3.3} we find that if we assume Conjecture \ref{conjecture2.3} then these algebraicity results hold for $\sum g_{i}(u_{i},v_{i})+x^{3}+y^{3}+xyz$ whenever $F$ is finite over $Z/2$. 

\end{theorem}

In Theorem \ref{theorem3.1} it is possible in theory, once the $r_{i,j}$, the $s_{i,j}$ and the $c_{j}$ are known, to get a polynomial relation between $w$ and $\mathcal{S}_{t+E_{1}}$ and compute $\mu(t+E_{1})$ by using the methods of \cite{3}. This is daunting in practice but we'll give one interesting partial result. Let $M$ be the smallest subspace of $X/(Q+Q\cdot t)$ that contains the image of $\beta_{1}$ and is stable under $T_{0}$ and $T_{1}$; our hypotheses show it to be finite dimensional. If $J_{0}$ and $J_{1}$ are maps $M\rightarrow M$ let $\Psi_{J_{0},J_{1}}(x,w)$ be the $2$-variable polynomial $\det |xI-w^{2}(J_{0}+xJ_{1})(J_{1}+xJ_{0})|$.

\begin{theorem}
\label{theorem3.5}
In the situation of Theorem \ref{theorem3.1}, $\sum E_{1}(2^{-n})(8w)^{n}$ lies in the splitting field over $Q(w)$ of $\Psi_{T_{0},T_{1}}(x,w)$.

\end{theorem}

\begin{proof}
We adopt the notation of Theorem \ref{theorem3.1} and its proof. $\sum E_{1}(2^{-n})(8w)^{n}=\sum l(T^{n}(E_{1}))w^{n}$, and Theorem 5.11 of \cite{3} shows that this power series lies in a certain extension $\mathcal{L}$ of $Q(w)$ constructed from the matrices $A$, $B$ and $C$. We saw in \cite{3} that $\mathcal{L} \subset $ a splitting field over $Q(w)$ of $\det |xI_{s}-w(Ax^{2}+Bx+C)|$. This last matrix is 
$$\begin{pmatrix}
xI_{l}&-wx(R+xS)\hspace{.5em}\\[-1.5ex]
\hspace{.5em}-w(S+xR)&xI_{l}
\end{pmatrix}.$$ So our determinant is just
$$x^{l}\det\begin{pmatrix}
I_{l}&-w(R+xS)\hspace{.5em}\\[-1.5ex]
\hspace{.5em} -w(S+xR)&xI_{l}
\end{pmatrix}.$$  Since $R$ and $S$ give the action of $T_{0}$ and $T_{1}$ on $M$, this last determinant is $\Psi_{T_{0},T_{1}}(x,w)$.
\qed
\end{proof}

\section{A (very) partially worked example}
\label{section4}

Suppose $\beta_{1}=\phi_{g}$ with $g=u^{6}+u^{3}v^{3}+v^{6}$. The methods of \cite{4} show that $M$ is five dimensional, that the action of $4T_{0}$ on $M$ is given by $\beta_{1}\rightarrow \beta_{2}\rightarrow \beta_{3}\rightarrow \beta_{1}$, $\beta_{4}\rightarrow \beta_{5}\rightarrow 0$, and that the action of $4T_{1}$ is given by $\beta_{5}\rightarrow \beta_{4}\rightarrow \beta_{3}\rightarrow \beta_{5}$, $\beta_{2}\rightarrow \beta_{1}\rightarrow 0$.  A Maple calculation then shows that $\Psi_{4T_{0},4T_{1}}(x,w)=-x^{2}\Psi^{*}$ where $\Psi^{*}$ is the reciprocal polynomial $w^{10}(x^{6}+1)-(2w^{8}+w^{4})(x^{5}+x)-(2w^{8}-3w^{6}-2w^{2})(x^{4}+x^{2})+(2w^{10}-w^{8}+2w^{6}-4w^{4}-1)x^{3}$. In an algebraic closure of $Q(w)$ let $\rho$, $\sigma$ and $\tau$ be the roots of $\Psi^{*}$ having positive ord; the other 3 roots are $\rho^{-1}$, $\sigma^{-1}$ and $\tau^{-1}$. The Galois group of $\Psi^{*}$ over $Q(w)$ has order 48 and consists of those permutations of the roots that permute the sets $\{\rho,\rho^{-1}\}$, $\{\sigma,\sigma^{-1}\}$, $\{\tau,\tau^{-1}\}$ among themselves.

Now Theorem \ref{theorem3.5} shows that $\sum E_{1}(2^{-n})(32w)^{n}$ is in a splitting field of $\Psi^{*}$ over $Q(w)$. But as we saw in \cite{3}, the field $\mathcal{L}$ attached to the matrices $A$, $B$ and $C$ sits inside a certain subfield of the splitting field of $\det |xI_{s}-w(Ax^{2}+Bx+C)|$. In our case $\mathcal{L}\subset$ the degree 8 extension of $Q(w)$ corresponding to the subgroup of the Galois group that stabilizes the set $\{\rho,\sigma,\tau\}$. Let $u_{1}=w^{10}(\rho-\rho^{-1})(\sigma-\sigma^{-1})(\tau-\tau^{-1})$ and $u_{2}=w^{10}(\rho\sigma\tau + \rho^{-1}\sigma^{-1}\tau^{-1})$. Using Galois theory we find that $u_{1}^{2}$ is in $Q(w)$, that $u_{2}$ has degree 4 over $Q(w)$, and that $u_{1}$ and $u_{2}$ generate the degree 8 extension of $Q(w)$ mentioned above.

So $\sum E_{1}(2^{-n})(32w)^{n}$ lies in $Q(w,u_{1},u_{2})$. A short calculation shows that $u_{1}^{2}=(w^{2}-1)^{2}(w^{2}+1)^{4}((1-w^{2})^{2}-4w^{6})$. One can also write down an irreducible equation for $u_{2}$ over $Q(w)$ but it's messy. (Some of the primes of $Q[w]$ that ramify in $Q(w,u_{2})$ are $(1-w^{2}+2w^{3})$, $(1-w^{2}-2w^{3})$ and $(4+8w^{2}-4w^{4}-12w^{6}-23w^{8}-18w^{10}+81w^{12}+108w^{14})$). Now the only fields between $Q(w)$ and $Q(w,u_{1},u_{2})$ are $Q(w)$, $Q(w,u_{1})=Q(w,\sqrt{(1-w^{2})^{2}-4w^{6}})$, $Q(w,u_{2})$ and $Q(w,u_{1},u_{2})$. So $\sum E_{1}(2^{-n})(32w)^{n}$, and consequently the conjectured Hilbert-Kunz series of $u^{6}+u^{3}v^{3}+v^{6}+x^{3}+y^{3}+xyz$, generates one of these 4 extensions of $Q(w)$. I think it generates the full degree 8 extension, but verifying this would be a very nasty computation.

Now consider the integral closure of $Q[w]$ in $Q(w,u_{1},u_{2})$. There is just one prime ideal in this ring lying over $(1-16w)$, and the argument of Theorem \ref{theorem1.14} shows that $\mu(t+E_{1})$, the putative Hilbert-Kunz multiplicity of $u^{6}+u^{3}v^{3}+v^{6}+x^{3}+y^{3}+xyz$ lies in the residue class field of this ideal.

The residue-class field is a degree 8 extension of $Q$ generated by the images, $\bar{u}_{1}$ and $\bar{u}_{2}$ of $u_{1}$ and $u_{2}$. $Q(\bar{u}_{1})$ is just $Q(\sqrt{(13)(157)(2039)})$, while $Q(\bar{u}_{2})$ is a degree 4 extension of $Q$ with discriminant $2^{2}\cdot 3^{3}\cdot 5^{2}\cdot 13^{2}\cdot 17^{2}\cdot 31\cdot 157^{2}\cdot 2039^{2}\cdot 780854102129687$. The only subfields of $Q(\bar{u}_{1},\bar{u}_{2})$ are $Q$, $Q(\bar{u}_{1})$, $Q(\bar{u}_{2})$ and $Q(\bar{u}_{1},\bar{u}_{2})$. So $\mu(t+E_{1})$ generates one of these 4 extensions of $Q$. My belief is that it generates the full degree 8 extension.



\label{}




\begin{thebibliography}{00}




\bibitem{1} C. Han, P. Monsky, Some surprising Hilbert-Kunz functions, Math. Z. 214 (1993), 119--135.

\bibitem{2} P. Monsky, Rationality of Hilbert-Kunz multiplicities: a likely counterexample, Michigan Math.\ J. 57 (2008), 605--613.

\bibitem{3} P. Monsky, Generating functions attached to some infinite matrices, Preprint (2009), arXiv:math.CO/0906.1836.

\bibitem{4} P. Monsky, P. Teixeira, $p$-Fractals and power series I, J. Algebra 280 (2004), 505--536.

\bibitem{5} P. Monsky, P. Teixeira, $p$-Fractals and power series II, J. Algebra 304 (2006), 237--255.

\end{thebibliography}
\end{document}